\documentclass[final]{iise}

\usepackage{multirow}
\usepackage{longtable}
\usepackage{booktabs}
\usepackage{soul}
\allowdisplaybreaks

\conference{Proceedings of the IISE Annual Conference \& Expo 2024 \\
A. Brown Greer, C. Contardo, J.-M. Frayret, eds.}

\title{\titlesize Multi-Period Stochastic Logistic Hub Capacity Planning\\ for Relay Transportation}

\author{
Xiaoyue Liu, Jingze Li, Mathieu Dahan, Benoit Montreuil\\Physical Internet Center, H. Milton Stewart School of Industrial \& Systems Engineering,\\Georgia Institute of Technology, Atlanta, GA 30332\\
}
\authorlist{Xiaoyue Liu, Jingze Li, Mathieu Dahan, Benoit Montreuil}
\abstractID{7118}

\begin{document}
\maketitle

\begin{abstract}
{\small This study focuses on relay transport carriers (RTCs) that contract with hub providers to lease hub capacity and employ relay transportation via hubs. It enables long-haul freight shipments to be transported by multiple short-haul drivers commuting between fixed-base hubs, promoting a driver-friendly approach. Inspired by Physical Internet, our paper addresses the multi-period capacity planning of logistic hubs within relay networks, accounting for uncertainty in demand and travel times. We model the problem as a two-stage stochastic optimization to determine the dynamic logistic hub throughput capacities for each planning period, ensuring the fulfillment of logistic demand while simultaneously minimizing both hub and transportation costs. This optimization problem falls within the NP-hard complexity class. To alleviate the inherent challenges in solving this problem, we employ a scenario reduction algorithm based on the fast forward selection (FFS) method to reduce computational effort while preserving approximation quality. Experiments with an automotive-delivery RTC in the Southeastern US demonstrate that our capacity planning model enables RTCs to proactively respond to dynamic circumstances, curtail avoidable expenditures, and enhance overall logistical efficiency. 
}
\end{abstract}

\section*{Keywords}
Multi-Period Stochastic Capacity Planning, Relay Transportation, Hub Capacity-Routing Problem, Physical Internet

\section{Introduction}

As reported in \cite{ata2021}, truck driver shortage reached a historic high of over 80,000 in 2021, and predictions indicate it may exceed 160,000 by 2030. The traditional end-to-end transportation heavily depends on long-haul freight movement. Covering extensive distances in a single trip requires drivers to endure prolonged periods away from home, which negatively impacts the physical and mental health of drivers, constituting a major factor in truck driver shortage \citep{li2022trucker,li2023multi}. To address these problems, many researchers are exploring a shift from conventional end-to-end transportation to relay transportation \citep{montreuil2011toward, kulkarni2022resilient}. The key idea behind relay transportation is to position logistic hubs along transport routes. Then the tasks of long-haul freight shipments can be distributed among multiple short-haul drivers commuting between fixed-base hubs. This approach enables truckers to operate near their residences and return home everyday. 

However, one of the main challenges in adopting relay transportation is the need for a large number of physical facilities - logistic hubs, which is especially unfriendly to small and medium-sized carriers. The introduction of the Physical Internet (PI) offers a promising avenue for addressing this challenge. Proposed by [4], PI presents an innovative approach to shaping supply chains, guiding them towards sustainability, economic efficiency, and social responsibility. One of the core concepts of PI is to promote an open and cooperative sharing economy. Distinguishing from the ordinary sharing economy, the PI infrastructure is open to multiple parties rather than exclusively dedicated to a single company or group of companies. Building on this concept, \cite{furtado2014impact} utilize a multi-agent simulation to demonstrate that resource-sharing transportation within the framework of PI enhances the performance of transport operations across financial, operational, social, and environmental dimensions.

Motivated by the innovative concept of PI, we examine a setting where an RTC aims to contract with hub providers to secure space for relay transportation under demand variations and uncertain travel time caused by network disruptions. Within this context, we define a multi-period stochastic hub capacity planning with routing problem (referred to as MSHCRP) to determine dynamic hub throughput capacities as well as shipment routing at each planning period. We then formulate a two-stage stochastic optimization model for the MSHCRP, aiming to minimize overall expenses while addressing both demand and travel time uncertainty. However, the resulting mixed-integer programming (MIP) poses computational challenges for real-world problems. To tackle this issue, we adopt a scenario reduction algorithm proposed by \cite{heitsch2003scenario} to approximate underlying continuous distributions for demand and travel time.

The rest of this paper is organized as follows. Section \ref{Sec:rr} briefly summarizes literature reviews of the relevant works. In Section \ref{Sec:pd}, we describe the problem setting and state research assumptions. Section \ref{Sec:model} formulates a two-stage stochastic program model and presents the scenario generation and reduction method to alleviate the computational challenge. Section \ref{Sec:er} delineates the dataset and experiments, extracting insights from the obtained results. Finally, Section \ref{Sec:concl} concludes this paper and discusses future research avenues.

\section{Related Research}
\label{Sec:rr}
The hub location-routing problem (HLRP) has received considerable attention in academic research for decades. Introduced by \cite{nagy1998many}, HLRP comprises two levels of decisions: a strategic level focused on determining hub locations and an operational level dealing with route planning. However, a notable critique of classical (single-period) HLRP is the disparity between the long-term nature of strategic decisions and the short-term nature of operational decisions \citep{alumur2021perspectives}. To bridge this temporal gap, the multi-period HLRP emerges as a viable solution, offering more frequent planning intervals for routing decisions. Although this multi-period approach is not novel, relevant research within the HLRP literature remains limited. Previous work by \cite{aloullal2023multi} explored scenarios involving progressive hub openings or closures, and studies by \cite{khaleghi2023multi} considered dynamic modular hub capacity that can gradually expand over the planning horizon. These considerations reflect real-world constraints, where budgetary limitations and resource availability may prohibit the simultaneous setup of an entire hub network. Another issue identified in the existing literature on HLRP is the heavy reliance on long-haul deliveries due to the significant costs associated with constructing hubs. This results in poor working conditions for drivers and further exacerbates the ongoing truck driver shortage. To alleviate this problem, we suggest an innovative solution where carriers can dynamically lease hub capacity through short-term contracts and use these hubs as simple relay facilities, which are not only cost-effective but also easily adjustable. This variant of multi-period HLRP allows carriers to manage long-haul shipments efficiently using short-haul drivers, promoting a driver-friendly approach. Additionally, it provides greater flexibility for carriers to adapt to dynamic environments.

The increasing complexities of real-world scenarios highlight the need to incorporate stochastic elements into problem modeling. In the literature on HLRP, two main types of uncertainties are usually considered separately: demand uncertainty, reflecting changes or fuzziness in demand, and travel time uncertainty, associated with factors such as traffic congestion or network failures. For instance, \cite{khaleghi2023multi} utilized a polynomial function to capture time-dependent demand dynamics, addressing a multi-period HLRP through various meta-heuristic algorithms. \cite{yang2009stochastic} addressed seasonal fluctuations in the air freight market, formulating air hub location and flight routing as a two-stage stochastic program model. In a different context, \cite{butun2021capacitated} explored HLRP with congestion effects in cargo delivery, while \cite{chaharsooghi2017adaptive} proposed an HLRP model based on two-stage stochastic optimization under the risk of hub disruptions. However, to the best of our knowledge, limited effort has been made to integrate both demand and travel time uncertainty within a multi-period framework, let alone adapting it to the context of relay-based transportation. Our study aims to bridge this gap in the HLRP literature by formulating a two-stage stochastic optimization model that simultaneously addresses two distinct uncertainties in a multi-period and relay-based setting. In the experimental phase, we conduct numerical experiments based on real automotive data, focusing on critical uncertainty parameters. These experiments yield valuable insights into the impact of demand and travel time uncertainty on network robustness, cost-effectiveness, delivery efficiency, and environmental impacts.

\section{Problem Description}
\label{Sec:pd}
This work tackles a logistic hub capacity planning challenge within a relay-based transportation network under a multi-period and multi-uncertainty environment. The key objective is to efficiently serve logistical demands while minimizing hub costs and transportation costs under uncertain conditions. In our study, logistic hubs primarily serve as facilities offering parking spots for driver waiting and resting, along with switch zones and buffer spots for driver and trailer swaps to support relay transportation. This definition aligns with the design of PI transit centers \citep{meller2013functional}. The contracted capacity in this scenario denotes the throughput capacity of logistic hubs, involving modular factors such as the number of parking spots and trailer swap docks. Specifically, we assume that the RTC makes one-year contracts with hub providers, indicating a total planning horizon of one year. Within this time frame, we explore a case where carriers have the flexibility to modify their contracted throughput capacities across different planning seasons in response to evolving demand needs, albeit with associated capacity-changing penalties. In Section \ref{Sec:er}, we conduct experiments to delve into the advantages of this dynamic (i.e., multi-period) setting. 

\begin{figure}[h!]
    \centering
    \includegraphics[width=13.5cm]{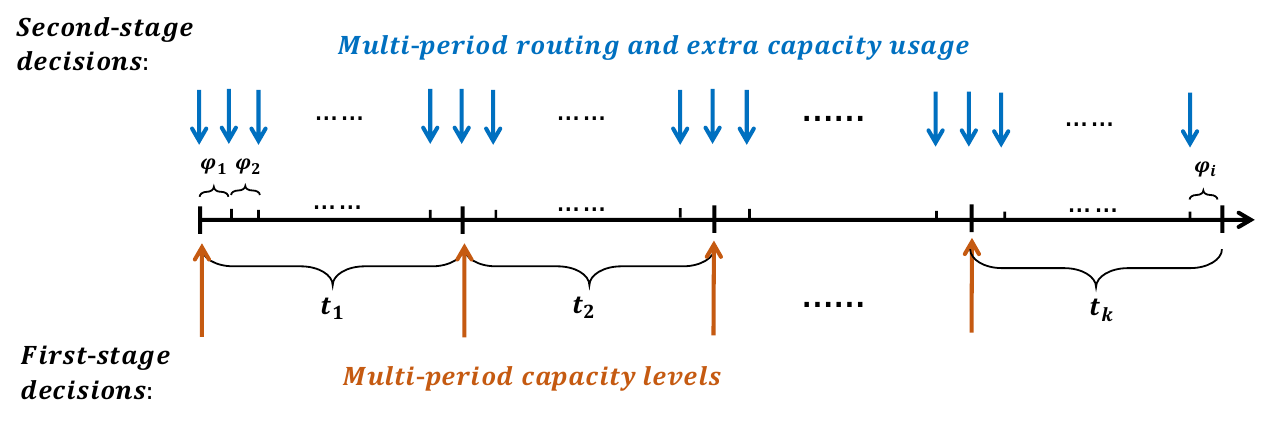}
    \caption{Timeline of the multi-period stochastic capacity planning model}
    \label{fig:timeline}
\end{figure}

A primary challenge lies in accounting for different time spans for various decision types. To address this, we formulate the problem using a two-stage stochastic program that divides decisions into two stages: The first stage involves logistic hub capacity level decisions, while the second stage encompasses freight routing decisions and extra capacity decisions. As depicted in Figure \ref{fig:timeline}, capacity level decisions are made at the beginning of each planning season $t_k$, whereas routing and extra capacity usage decisions are made for each operational period $\varphi_i$ after realizing demand and geographical disruptions. We assume all shipments are dispatched at the beginning of each operational period and are fully routed through the network, without occupying the transportation resources of the subsequent operational period.

Another major challenge is to determine logistic hub capacity in the face of uncertainty. When signing a capacity contract, the demand for future periods remains unknown. However, it is the demand that primarily influences capacity decisions. To model demand uncertainty, we assume that the carrier has historical data reflecting demand distribution from past service records, and model this uncertainty through probabilistic scenarios. Yet, even with sufficient data for scenario generation, incorporating too many scenarios can pose a computational burden. To balance this trade-off, we implement a scenario reduction to reduce computational effort while preserving the quality of the approximation. 

\section{Problem Formulation}
\label{Sec:model}
We model this problem as a two-stage stochastic program to strategically determine logistic hub throughput capacities for each planning season, ensuring the fulfillment of logistic needs while minimizing both hub costs and transportation costs. The first stage decisions involve determining the throughput capacity levels of each logistic hub. In the second stage, following the realization of demand and geographical disruptions, decisions are made regarding freight shipment routing within the network and the utilization of extra throughput capacity at each logistic hub. The model notations are displayed in Table \ref{tab:notations} and its formulation are as follows:
{
\small
\begin{longtable}[H]{llll}
    \caption{Notation for the multi-period stochastic capacity planning model}\\ 
    \vspace{-0.7cm}\\
    \hline
    Sets & &Parameters &\\
    \hline
    $\boldsymbol{O}$ & Set of origins, indexed by $o$ & $q_{o d}(\varphi,\omega)$ & Demand from origin $o$ to destination $d$ in \\
    $\boldsymbol{D}$ & Set of destinations, indexed by $d$ && operational period $\varphi$ under scenario $\omega$\\
    $\boldsymbol{H}$ & Set of logistic hubs, indexed by $h$ & $t_{a}(\varphi, \omega)$ & Travel time on arc $a$ in operational period \\
    $\boldsymbol{P}$ & Set of demand OD pairs, indexed by $(o,d)$ && $\varphi$ under scenario $\omega$ (affected by disruptions)\\
    $\boldsymbol{A}$ & Set of arcs, indexed by $a$ or $(i,j)$ & $c^H_{h l_1 l_2}$ & Hub cost to change capacity level from $l_1$ to \\
    $\boldsymbol{T}$ & Set of planning seasons, indexed by $t$ && $l_2$ and operate hub $h$ at capacity level $l_2$\\
    $\boldsymbol{T_t}$ & Set of operational periods in season $t$, indexed by $\varphi$ & $c^F$ & Flow unit cost per hour\\
    $\boldsymbol{W_t}$ & Set of scenarios in season $t$, indexed by $\omega$ & $c^P_h$ & Penalty cost for using extra capacity at hub $h$\\
    $\boldsymbol{L}$ & Set of capacity levels, indexed by $l$ & $u_{hl}$ & Hub capacity at capacity level $l$ at hub $h$\\
    $\boldsymbol{\delta^{+}(i)}$ & Set of arcs exiting node $i$ & $m$ & Number of days during an operational period \\
    $\boldsymbol{\delta^{-}(i)}$ & Set of arcs entering node $i$ & $p_\omega$ & Probability of scenario $\omega$ \\
    \hline
    \label{tab:notations}
\end{longtable}
}
\vspace*{-32pt}
{
\small
\begin{longtable}[H]{ll}
    Variables\\
    \hline
    $X_{h l_1 l_2}(t)$ & Whether to change capacity level  from $l_1$ to $l_2$ at hub $h$ at the beginning of planning season $t$ \qquad \qquad \qquad \quad\\
    $F_{oda}(\varphi,\omega)$ & Freight volume from origin $o$ to destination $d$ on arc $a$ in period $\varphi$ under scenario $\omega$\\
    $E_{h}(\varphi,\omega)$ & Extra capacity usage at hub $h$ in operational period $\varphi$ under scenario $\omega$\\
    \hline
    \addtocounter{table}{-1}
\end{longtable}
}
\vspace*{-24pt}
\begin{flalign}
    & \min \sum_{t \in \boldsymbol{T}}\sum_{h \in \boldsymbol{H}} \sum_{l_1 \in \boldsymbol{L}} \sum_{l_2 \in \boldsymbol{L}} c_{h l_1 l_2}^H \cdot X_{h l_1 l_2}(t) + \sum_{t \in \boldsymbol{T}} \sum_{\omega \in \boldsymbol{W_t}} \sum_{\varphi \in \boldsymbol{T_t}} p_\omega\left(\sum_{a \in \boldsymbol{A}} \sum_{(o, d) \in \boldsymbol{P}} c^F \cdot t_a(\varphi, \omega) \cdot F_{o d a}(\varphi, \omega)+\sum_{h \in \boldsymbol{H}} c_h^P \cdot E_h(\varphi, \omega)\right) \label{eq:1} \\ 
    & s.t. \notag
\end{flalign}
\vspace*{-15pt}
\begin{flalign}
    & \sum_{a \in \boldsymbol{\delta^{+}(o)}} F_{o d a}(\varphi, \omega)=q_{o d}(\varphi, \omega), \quad \forall(o, d) \in \boldsymbol{P}, \forall \varphi \in \boldsymbol{T}_{\boldsymbol{t}}, \forall \omega \in \boldsymbol{W}_{\boldsymbol{t}}, \forall t \in \boldsymbol{T} \label{eq:2}\\
    & \sum_{a \in \boldsymbol{\delta^{-}(d)}} F_{o d a}(\varphi, \omega)=q_{o d}(\varphi, \omega), \quad \forall(o, d) \in \boldsymbol{P}, \forall \varphi \in \boldsymbol{T}_{\boldsymbol{t}}, \forall \omega \in \boldsymbol{W}_t, \forall t \in \boldsymbol{T} \label{eq:3}\\
    & \sum_{a \in \boldsymbol{\delta^{+}(h)}} F_{o d a}(\varphi, \omega)=\sum_{a \in \boldsymbol{\delta^{-}(h)}} F_{o d a}(\varphi, \omega), \quad \forall h \in \boldsymbol{H}, \forall(o, d) \in \boldsymbol{P}, \forall \varphi \in \boldsymbol{T}_{\boldsymbol{t}}, \forall \omega \in \boldsymbol{W}_t, \forall t \in \boldsymbol{T} \label{eq:4}\\
    & \sum_{(o, d) \in \boldsymbol{P}} \sum_{a \in \boldsymbol{\delta^{-}(h)}} F_{o d a}(\varphi, \omega) \leq m \cdot\left(\sum_{l_1 \in \boldsymbol{L}} \sum_{l_2 \in \boldsymbol{L}} u_{h l_2} \cdot X_{h l_1 l_2}(t)+E_h(\varphi, \omega)\right), \quad \forall h \in \boldsymbol{H}, \forall \varphi \in \boldsymbol{T}_{\boldsymbol{t}}, \forall \omega \in \boldsymbol{W}_{\boldsymbol{t}}, \forall t \in \boldsymbol{T} \label{eq:5}\\
    & \sum_{l \in \boldsymbol{L}} X_{h 0 l}(0)=1, \quad \forall h \in \boldsymbol{H} \label{eq:6}\\
    & \sum_{l_1 \in \boldsymbol{L}} \sum_{l_2 \in \boldsymbol{L}} X_{h l_1 l_2}(t)=1, \quad \forall h \in \boldsymbol{H}, \forall t \in \boldsymbol{T} \label{eq:7}\\
    & \sum_{l_1 \in \boldsymbol{L}} X_{h l_1 l}(t-1)=\sum_{l_2 \in \boldsymbol{L}} X_{h l} l_2(t), \quad \forall l \in \boldsymbol{L}, \forall h \in \boldsymbol{H}, \forall t \in \boldsymbol{T} \backslash\{\mathbf{0}\} \label{eq:8}
\end{flalign}
\vspace*{-2pt}
\begin{flalign}
    & X_{h l_1 l_2}(t) \in\{0,1\}, E_h(\varphi, \omega) \geq 0, F_{o d a}(\varphi, \omega) \geq 0, \notag\\
    & \hspace{179.9 pt minus 1fil}   \forall a \in \boldsymbol{A}, \forall(o, d) \in \boldsymbol{P}, \forall l_1, l_2 \in \boldsymbol{L}, \forall h \in \boldsymbol{H}, \forall \varphi \in \boldsymbol{T}_{\boldsymbol{t}}, \forall \omega \in \boldsymbol{W}_{\boldsymbol{t}}, \forall t \in \boldsymbol{T} \label{eq:9}
\end{flalign}

The objective function (\ref{eq:1}) is designed to minimize the expected total costs over the scenario set during the entire planning horizon. Flow conservation constraints (\ref{eq:2})-(\ref{eq:4}) ensure the proper routing of each freight shipment from its origin to its destination. Constraints (\ref{eq:5}) prevent capacity restriction violations for each logistic hub. Constraints (\ref{eq:6}) affirm that the initial throughput capacity at each logistic hub is empty. Constraints (\ref{eq:7}) and (\ref{eq:8}) assure the selection of exactly one hub capacity deployment strategy in each planning season, with consecutive changes in capacity levels for each hub. Lastly, constraints (\ref{eq:9}) specify the domains of the model variables.

In this study, we tackle two types of uncertainty: demand and travel time uncertainty. To address computational challenges in stochastic models with continuous distributions, we adopt a common approach of modeling uncertainty using discrete probabilistic scenarios. Specifically, we utilize a Poisson distribution to generate weekly demand from the Freight Analysis Framework (FAF) yearly demand dataset \cite{faf}. The Poisson arrival rate is determined by the historical monthly demand share, evenly distributed across four weeks within each month. This process is iterated 50 times to create 50 scenarios, all sharing the same arrival rate. For travel time uncertainty, we consider two key parameters: disruption rate and intensity. We assume that hubs operate independently and have a probability of being disrupted based on a given rate. For example, a disruption rate of 0.05 and intensity of 1.5 mean that each hub faces a 5\% probability of disruption, and the travel time on the arc connecting to the interruption hub is 1.5 times the normal travel time. This simulates scenarios such as construction and maintenance within a hub, occurring independently.

The problem size, in terms of variables and constraints, significantly expands with an increasing number of scenarios. Given our focus on a large-scale relay network, the computational load becomes substantial with a large number of scenarios. To mitigate this issue, a common approach involves approximating the original set of scenarios with a smaller representative subset. In this work, we employ the fast forward selection (FFS) method, which is a widely used scenario reduction approach introduced by \cite{heitsch2003scenario}. This FFS method selects the 10 scenarios out of the original scenario set with the shortest Kantorovich distance from the unselected set, making it the most representative of the unselected group, thus reducing computational burden while preserving solution quality.

\section{Experimental Results}
\label{Sec:er}

\begin{figure}[h!]
    \centering
    \includegraphics[width=15.5cm]{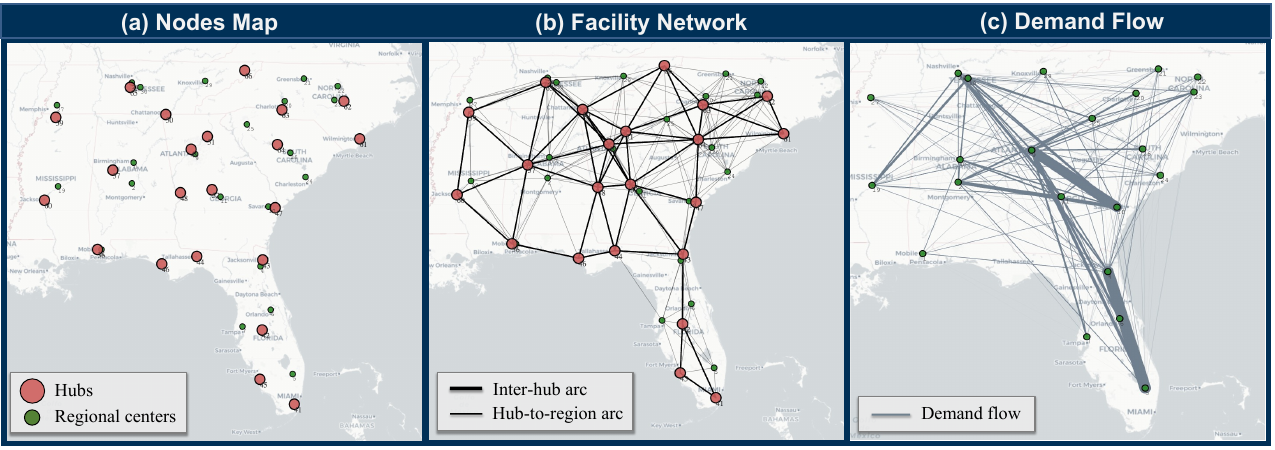}
    \caption{Nodal map (a), facility network (b) and demand flow (c) for multiple automotive OEMs in Southeast USA}
    \label{fig:dataset_region2}
\end{figure}

In this section, we employ a yearly demand dataset for finished automotive vehicles sourced from the FAF database, with a specific focus on the southeast US. We conduct experiments focusing on critical uncertainty parameters, aiming to gain a comprehensive understanding of their influence on the economic and social performance of the system. The network nodes, as depicted in Figure \ref{fig:dataset_region2}(a), comprise 22 logistic hubs and 22 region centers. Figure \ref{fig:dataset_region2}(b) provides an overview of the relay network structure, with edges configured to ensure travel time along each route is limited to 5.5 hours. This design allows short-haul drivers to complete a round trip within 11 hours, in compliance with the daily driving hour limit set by US regulations. Lastly, Figure \ref{fig:dataset_region2}(c) illustrates the yearly demand flow across the southeast US, reflecting the distribution of freight volume. 

To assess the impact of demand and travel time uncertainty on logistic hub capacity planning, we run the MSHCRP model across three uncertainty levels: low, intermediate, and high. Each level corresponds to a specific scenario set, generated using the method outlined in Section \ref{Sec:model}. Specifically, at the low uncertainty level, there is no disruption in the network, focusing exclusively on demand uncertainty. Moving to the intermediate uncertainty level, we introduce disruption uncertainty, where disruption scenarios are generated with a disruption rate of 0.05 and a disruption intensity of 1.5. Lastly, at the high uncertainty level, we examine a scenario with a disruption rate of 0.05 and a disruption intensity of 15, representing a significant delay penalty cost. To gain deeper insights into the benefits of dynamic modeling, we further compare the MSHCRP model with its static counterpart, where contracted capacities, once determined by the model, remain fixed throughout the entire planning horizon. All results are summarized in Table \ref{tab:results}. 

\begin{table}[h!]
    \centering
    \caption{Performance Metrics under Different Uncertainty Levels}
    \vspace{-0.3cm}
    \resizebox{\linewidth}{!}{
    \begin{tabular}{lcccccccc}
        \toprule
        Uncertainty Levels & Low &&& Intermediate &&& High\\ \cmidrule{2-3} \cmidrule{5-6} \cmidrule{8-9}
        Model types & Static & Dynamic && Static & Dynamic && Static & Dynamic \\
        \midrule
        Total Contracted Throughput Capacity & 4,080 & 3,720 && 4,200 & 3,870 && 4,320 & 3,870 \\
        Average Hub Throughput Capacity & 72.86 & 66.43 && 70.00 & 64.50 && 72.00 & 63.80 \\
        Average Number of Hubs & 14.00 & 14.00 && 15.00 & 15.00 && 15.00 & 15.00 \\
        Average Hub Network Connectivity & 4.43 & 4.43 && 4.67 & 4.67 && 4.67 & 4.52 \\
        Hub Costs (\$) & 11,166,155 & 9,759,231 && 11,281,953 & 9,953,695 && 12,108,207 & 10,885,179 \\
        Extra Capacity Usage Percentage (\%) & 15.19 & 5.67 && 13.78 & 5.27 && 19.04 & 11.54 \\
        Extra Capacity Usage Frequency (\%) & 3.53 & 1.65 && 3.18 & 1.39 && 4.12 & 3.48 \\
        Transportation Costs (\$) & 34,558,081 & 34,623,956 && 35,185,264 & 35,249,483 && 36,004,177 & 36,138,348 \\
        Disruption Time Percentage (\%) & 0.00 & 0.00 && 4.77 & 4.92 && 4.70 & 5.31 \\
        Disrupted Edge Usage Frequency (\%) & 0.00 & 0.00 && 3.12 & 3.24 && 1.23 & 1.40 \\
        Fuel Consumption (tons) & 38,878 & 38,952 && 39,583 & 39,656 && 40,505 & 40,656 \\
        CO2 Emissions (metric tons) & 107,924 & 108,129 && 109,882 & 110,083 && 112,440 & 112,859 \\
        Total Costs (\$) & 45,724,236 & 44,383,187 && 46,467,217 & 45,203,177 && 48,112,384 & 47,023,527 \\
        \bottomrule
    \end{tabular}
    \label{tab:results}
    }
\end{table}

Firstly, we focus on the dynamic model (specifically the MSHCRP model) under the changes in uncertainty levels. As shown in Table \ref{tab:results}, the total costs exhibit an increase from \$44,383,187 to \$47,023,527 with the escalation of uncertainty levels. Simultaneously, we can observe a decline in the average hub throughput capacity and a rise in the number of hubs. This suggests a trend within the model towards opening more hubs with a smaller design in the face of larger uncertainty, strategically mitigating hub disruption risks. In the transportation dimension, the fuel consumption increases from 108,129 tons to 112,859 tons, indicating the model is inclined to prioritize alternative longer routes to avoid disrupted arcs, thereby minimizing their impacts on operations. 

Subsequently, we shift the focus to the comparison between the static model and dynamic model, as outlined in Table \ref{tab:results}. The findings reveal that, through dynamic capacity planning, the carrier can realize savings exceeding \$1,000,000 in hub costs by adjusting their logistic hub capacities across various planning seasons in response to evolving operational needs. Despite the reduced contracted capacity, the dynamic model uses approximately one-half to one-third of extra capacity compared to the static model, showcasing robustness in dynamic environments. On the other hand, the reduction in total hub capacity within the dynamic model has a negligible impact on hub network connectivity, transportation metrics, and environmental metrics, all of which remain relatively consistent with the static model. Overall, the results demonstrate that adopting the dynamic setting enables RTCs to proactively respond to dynamic circumstances, curtail avoidable expenditures, and enhance comprehensive logistical efficiency.

\section{Conclusion}
\label{Sec:concl}
This study introduces a two-stage stochastic optimization model to tackle the multi-period stochastic hub capacity planning with routing problem (MSHCRP) in the context of relay transportation, which few researchers have paid attention to. Additionally, we leverage a scenario reduction algorithm based on FFS to approximate continuous distributions for demand and travel time. In the experimental section, we analyze the impact of uncertainties across four dimensions (network structure, hub contract, transportation, and environment) and compare the model with the static capacity planning model. The results demonstrate that dynamic capacity planning outperforms static capacity planning in terms of network robustness and cost-effectiveness, while maintaining strong performance in delivery efficiency and environmental impact.

Several promising avenues for future research can be identified. Firstly, investigating advanced computational methods like accelerated Benders decomposition and branch-and-price could enhance the solution approach, making it more efficient for large-scale instances. Secondly, the integration of relay-based transportation with emerging technologies, such as automatic and connected trucks, offers an exciting prospect for exploration. Thirdly, a more detailed categorization of hub capacity into distinct segments, such as loading zone capacity, parking zone capacity, and driver-switching zone capacity, could be beneficial. By implementing this approach, one can gain deeper insights, optimize resource allocation, and develop more effective strategies tailored to specific operational requirements.

\end{document}